\documentclass[10pt,twocolumn,twoside]{IEEEtran}

\usepackage{color}
\usepackage{amsfonts}
\usepackage{bbm}
\usepackage{mathrsfs}
\usepackage{amsmath}
\usepackage{graphics}
\usepackage{epsfig}
\usepackage{amssymb}
\usepackage{bm}
\usepackage{cite}
\usepackage{epstopdf}

\def\eq#1{\begin{equation}#1\end{equation}}
\def\rep#1{(\ref{#1})}
\newtheorem{theorem}{Theorem}

\newtheorem{lemma}{Lemma}
\newtheorem{remark}{Remark}

\newcommand{\R}{{\rm I\!R}}

\DeclareMathOperator{\diag}{diag}

\begin{document}

\title{\LARGE \bf Analysis of a Nonlinear Opinion Dynamics Model \\
with Biased Assimilation}
 \author{Weiguo Xia, Mengbin Ye, Ji Liu, Ming Cao, and Xi-Ming Sun
 \thanks{
 Weiguo Xia and Xi-Ming Sun are with the Key Laboratory of Intelligent Control and Optimization
for Industrial Equipment of Ministry of Education and School of Control Science and Engineering, Dalian University of Technology, China ({\tt \{wgxiaseu,sunxm\}@dlut.edu.cn}). Mengbin Ye is with the Optus-Curtin Centre of Excellence in Artificial Intelligence, Curtin University, Australia ({\tt mengbin.ye@curtin.edu.au}). Ji Liu is with the Department of Electrical and Computer Engineering, Stony Brook University, USA ({\tt ji.liu@stonybrook.edu}). Ming Cao is with the Faculty of Science and Engineering, ENTEG, University of Groningen, the Netherlands ({\tt m.cao@rug.nl}). }
 }

\maketitle

\begin{abstract}
This paper analyzes a nonlinear opinion dynamics model which generalizes the DeGroot model by introducing a bias parameter for each individual. The original DeGroot model is recovered when the bias parameter is equal to zero. The magnitude of this parameter reflects an individual's degree of bias when assimilating new opinions, and depending on the magnitude, an individual is said to have weak, intermediate, and strong bias. The opinions of the individuals lie between 0 and 1. It is shown that for strongly connected networks, the equilibria with all elements equal identically to the extreme value 0 or 1 is locally exponentially stable, while the equilibrium with all elements equal to the neutral consensus value of  $\frac{1}{2}$ is unstable. Regions of attraction for the extreme consensus equilibria are given. For the equilibrium consisting of both extreme values 0 and 1, which corresponds to opinion polarization according to the model, it is shown that the equilibrium is unstable for all strongly connected networks if individuals all have weak bias, becomes locally exponentially stable for complete and two-island networks if individuals all have strong bias, and its stability heavily depends on the network topology when individuals have intermediate bias. Analysis on star graphs and simulations show that additional equilibria may exist where individuals form clusters.
\end{abstract}


\section{Introduction}


There has been a persistent interest in theoretical sociology over the past decades in the modeling and study of opinion formation processes \cite{Aj01}. A variety of models have been proposed and studied to understand how the opinions of an interconnected social group evolve and how limiting phenomena arise, including consensus, polarization, and clustering.

{\color{black} The French--DeGroot model \cite{De74} is probably the most well-known (referred to as DeGroot henceforth for simplicity); each individual repeatedly updates his/her opinion to be a weighted average of the opinions of his/her neighbors (perhaps including him/herself), reflecting the subconscious human cognitive capability of taking convex combinations when processing new information \cite{An81}. The opinions of the individuals will eventually reach a consensus as long as the interaction network satisfies some appropriate connectivity requirements. Over the years, a number of discrete- and continuous-time variants of the DeGroot model have been proposed and studied extensively. Notable among them include the Friedkin--Johnsen model \cite{FrJo90}, the Hegselmann--Krause model \cite{HeKr02,BlHeTs09}, the Altafini model \cite{Al13,LiChBaBe17,XiCaJo16}, and the DeGroot-Friedkin model \cite{JiMiFrBu15,YeLiAnYuBa19}. For recent advances in the modelling of opinion dynamics on influence networks, see \cite{PrRo17,AnYe19}.


The phenomenon of extremization, which refers to the tendency for a group of individuals to eventually reach opinions that are more extreme than their initial inclination (and perhaps polarizing into two opposite camps), has become of increasing relevance in the modern age, and is the focus of research from several scientific communities \cite{Po14}. Most of the models discussed above are not able to predict polarization or extremization. Some models do exist in the literature which can predict polarization or extremization \cite{MaFlKi14,Du17,ShPrJoBaJo17}, but typically attributes this phenomenon to antagonistic interactions that increase in strength as the difference in opinions between individuals grow, and only \cite{ShPrJoBaJo17} has provided analysis for its proposed model.

{\color{black} Biased assimilation is the phenomenon in social psychology in which individuals tend to process new information with a bias towards their current position, accepting confirming evidence while evaluating disconfirming evidence critically \cite{LoRoLe79}. This can result in an individual developing a more extreme opinion when exposed to information from a confirming and disconfirming source \cite{MuDiLoFaGrPe02,TaLo06}.} A new generalization of the DeGroot model was recently proposed in \cite{DaGoLe13}, where a bias parameter helps to capture the cognitive processes described in the preceding two sentences. For homophilous networks, it has been shown that under some specific conditions, polarization arises if the individuals are sufficiently biased and consensus is reached under some other specific initial opinions for a small bias parameter close to zero \cite{DaGoLe13}. However, the situation that the system converges to consensus is rarely observed for other initial states.

{\color{black}In this paper, we focus on strongly connected networks and further examine the model proposed in \cite{DaGoLe13} to elucidate the role of biased assimilation in shaping opinion formation processes. The level of biased assimilation is captured by the scalar $b_i$ for individual $i$, and assumed to be heterogeneous among the individuals. First, we provide a detailed, quantitative argument to illustrate how biased assimilation is captured in the model when an individual is presented with two opposing opinions. As a consequence, we are able to clearly illustrate how the magnitude of $b_i$ determines whether individual $i$ has a weak, intermediate, or strong intensity of biased assimilation. We then concentrate on equilibria with meaningful social interpretations, such as \textit{extreme consensus} equilibria, the \textit{neutral consensus} equilibrium, and \textit{extreme polarization} equilibria. The role that $b_i$, viz. the intensity of the biased assimilation, plays in determining the (local) stability or instability of these equilibria is explicitly identified. For \textit{extreme consensus} equilibria, broad regions of attraction are obtained, and we identify further equilibria for star networks. This contrasts the work of \cite{DaGoLe13}, which identifies regions of convergence for the polarization equilibria but does not consider stability or other equilibria, and the work of \cite{ChQiLiQiBuSh19}, where stability results are only established for special classes of network topologies. Detailed discussions are provided for the findings on each of the above types of equilibria, and social interpretations and implications are examined.}

The rest of the paper is organized as follows. Section~\ref{model} introduces the biased opinion dynamics model. Section~\ref{result} analyzes the equilibria and their stability for the model under a general graph topology and discusses some classes of specific graphs with the proofs given in Section~\ref{analysis}. Section~\ref{simulation} provides several simulations to illustrate a rich set of possible dynamic behaviors possible, including some not covered in the analysis. Section~\ref{conclusion} concludes the paper.

{\it Notation:} For a positive  integer $N$, let $\bm{1}_N$ and $\bm{0}_N$ denote the $N$-dimensional all-one vector and all-zero vector, respectively.  Let $I_{N\times N}$ and $O_{N\times N}$ denote the $N\times N$ identity matrix and zero matrix, respectively. We will use the terms ``individual'' and ``agent'' interchangeably.

\section{The Model For Opinion Dynamics With Bias Assimilation} \label{model}

Consider a group of $N$ agents labeled by $1$ to $N$. Each agent can receive information only from its neighbors. The neighbor relationships
among the $N$ agents are characterized by an $N$-node directed graph represented by $\mathbb G=(\mathcal V, \mathcal E),$ where $ \mathcal V=\{1,\dots,N\}$ is the node set and $\mathcal E$ is the edge set. The graph is associated with a weight matrix $W=(w_{ij})_{N\times N}$ where the self-weight  $w_{ii}\ge 0$, and  {\color{black}if $(j,i)\in\mathcal E$, then $w_{ij}>0$.  Let $\mathcal{N}_i=\{j\in\mathcal V:(j,i)\in\mathcal E\}$ be the set of neighbors of agent $i$, representing the other agents $j$ that have influence on  $i$.} Note that no self-loop is allowed in the graph $\mathbb G$ and therefore $i\not\in\mathcal{N}_i$ for $i=1,\dots,N$, but the self-weight $w_{ii}$ can be positive. {\color{black}A directed path from node $p_1$ to node $p_k$ is a sequence of edges of the form $(p_1, p_{2}), (p_2, p_3), ...,(p_{k-1}, p_k)$ where $p_i \in \mathcal{V}$ are distinct and $(p_{i}, p_{i+1}) \in \mathcal{E}$. A graph is strongly connected if there is a path from every node to every other node, which is equivalent to $W$ being irreducible.}

Each agent $i$ has a real-valued opinion $x_i(k)$,  on a given issue being discussed, which may change over time $k$. {\color{black}At every discrete time instant $k = 0, 1, \hdots$, each agent $i$ updates its opinion by setting}
\eq{x_i(k+1) = \frac{w_{ii}x_i(k) + (x_i(k))^{b_i} s_i(k)}
{w_{ii} \!+\! (x_i(k))^{b_i} s_i(k) \!+\! (1-x_i(k))^{b_i} (d_i\!-\!s_i(k))},
\label{general}}
where $d_i=\sum_{j\in\mathcal{N}_i}w_{ij}$, $s_i(k)=\sum_{j\in\mathcal{N}_i} w_{ij}x_j(k)$, and $w_{ij}, i,j \in \mathcal V$, are the elements in the weight matrix $W$ representing the influence weights. The bias of agent $i$ is captured by the parameter $b_i$, and is assumed to be nonnegative except for a special scenario considered in the sequel. Observe that on the right hand side of \eqref{general}, the numerator is nonnegative and the denominator is greater than or equal to the numerator for any $x_i(k) \in [0,1]$. Thus, $x_i(0) \in [0,1]$ for all $i\in\mathcal V$ guarantees that $x_i(k) \in [0,1]$ for all $k \geq 0$ and $i\in \mathcal{V}$. {\color{black}We assume from here on that $x_i(0)\in[0,1]$ for all $i\in\mathcal V$, and $0$ and $1$ represent the extreme opinions of opposing points of view on the given topic, respectively. By way of example, suppose the issue being discussed was the statement ``recreational marijuana should be legalized'', then $x_i = 0$ and $x_i = 1$ correspond to individual $i$ totally opposing and totally supporting the legalization of marijuana. Consequently, $x_i(k)$ can be regarded as agent $i$'s degree of support at time $k$ for the extreme opinion represented by 1, and so correspondingly $1-x_i(k)$  can be regarded as the degree of support for the extreme opinion represented by 0. The reader is referred to \cite[Section 2.2]{AnYe19} for further details.}

{\color{black}We now give an intuitive explanation on how the model \eqref{general} captures bias assimilation, and provide quantitative arguments in the next subsection. Readers may also refer to \cite{DaGoLe13}, which first proposed the model. One can consider $s_i(k) = \sum_{j\in\mathcal{N}_i} w_{ij}x_j(k)$ and $d_i-s_i(k) = \sum_{j\in\mathcal{N}_i} w_{ij}(1-x_j(k))$ to be the weighted average support for the position represented by $1$ and 0, respectively. When $b_i > 0$ and supposing for example that $x_i(k) > 0$, \eqref{general} indicates that individual $i$ applies a larger weight of $x_i(k)^{b_i}$ to $s_i(k)$, and a smaller weight of $(1-x_i(k))^{b_i}$ to $d_i-s_i(k)$. This represents the biased assimilation phenomenon \cite{LoRoLe79}, which explains that individuals may process new information with a bias, being more readily inclined to accept evidence confirming their existing views while evaluating disconfirming evidence critically, perhaps even rejecting it. We remark that when $b_i=0$ for all $i\in\mathcal{V}$, \rep{general} simplifies to the classical DeGroot model \cite{De74}.}



\subsection{Exploring the Bias Parameter's Effect} \label{explanation}

{\color{black}In this section, we look closely at the effect of the bias parameter $b_i>0$ on the dynamics in (\ref{general}) and show that when $b_i > 0$, each individual assimilates new information with a bias towards information supporting his or her current opinion, and the value of $b_i$ determines the level of bias. To do so, we construct a specific example to understand the opinion update of a single individual $i$ in the presence of equal information from both ends of the opinion spectrum. The example imposes some additional assumptions, which are not restrictive; the same conclusions can be drawn with other similar assumptions.

Suppose that $w_{ij} =1$ for all $i,j \in \mathcal{V}$, and that the neighbors of $i$ have opinions that yield $s_i=\sum_{j\in\mathcal N_i}x_j=d_i/2\triangleq s$, i.e. there is an equal influence from $i$'s neighbors on both ends of the opinion spectrum\footnote{E.g., individual $i$ has two neighbors, one having an opinion value of 1, and the other an opinion value of 0.}.} The update equation \eqref{general} of individual $i$ can be rewritten as $$x_i(k+1)=p(b_i, x_i(k))$$ where
$$p(b_i,x_i)\triangleq\frac{x_i+x_i^{b_i}s}{1+x_i^{b_i}s+(1-x_i)^{b_i}s}.$$
Evidently, the  DeGroot update equation of individual $i$ is $x_i(k+1) = p(0,x_i(k)) = (x_i(k)+s)/(1+2s)$.

We will show that $p(b_i,x_i)>p(0,x_i)$ for all $b_i>0$ if $x_i\in(0.5,1)$, and $p(b_i,x_i)<p(0,x_i)$ for all $b_i>0$ if $x_i\in(0,0.5)$.

First, observe that
\begin{eqnarray}\label{eq:explain2}
p(b_i,x_i)-p(0,x_i)&=&\frac{g(b_i,x_i)}{(1+2s)(1+x_i^{b_i}+(1-x_i)^{b_i})},
\end{eqnarray}
where
\begin{eqnarray*}
g(b_i,x_i)&\triangleq& sx_i^{b_i}+2sx_i+2s^2x_i^{b_i}-s-sx_i^{b+1}-s^2x_i^{b_i}\\
&&-sx_i(1-x_i)^{b_i}-s^2(1-x_i)^{b_i}.
\end{eqnarray*}
The derivative of $g(b_i,x_i)$ with respect to $x_i$ is
\begin{eqnarray*}
&&\frac{\partial{g(b_i,x_i)}}{\partial{x}}=b_isx_i^{b_i-1}+2s+b_is^2x_i^{b_i-1}-(b_i+1)sx_i^{b_i}\\
&&\\&&\quad-s(1-x_i)^{b_i}+b_isx_i(1-x_i)^{b_i-1}+b_is^2(1-x_i)^{b_i-1}.
\end{eqnarray*}
Note that for all $b_i>0$ and $x_i\in(0,1)$, there holds (i) $b_isx_i^{b_i-1}-b_isx_i^{b_i}>0$, (ii) $s(1-x_i^{b_i})>0$, and (iii) $s(1-(1-x_i)^{b_i})>0$.
It follows that
$$b_isx_i^{b_i-1}+2s-(b_i+1)sx_i^{b_i}-s(1-x_i)^{b_i}>0,$$
and therefore $\frac{\partial{g(b_i,x_i)}}{\partial{x_i}}>0$. Combined with the fact that $g(b_i,\frac{1}{2})=0$ for $b_i > 0$, one has that $g(b_i,x_i)>0$ for all $b_i>0$ and $x_i\in(0.5,1)$, and $g(b_i,x_i)<0$ for all $b_i>0$ and $x_i\in(0,0.5)$. It follows from (\ref{eq:explain2}) that $p(b_i,x_i)>p(0,x_i)$ for all $b_i>0$ and $x_i\in(0.5,1)$, and $p(b_i,x_i)<p(0,x_i)$ for all $b_i>0$ and $x_i\in(0,0.5)$.

Note that when $b_i=1$, one has $p(1,x_i)=x_i$. When $b_i>1$, we obtain $p(b_i,x_i)>x_i$
for $x_i\in(0,5,1)$ and $p(b_i,x_i)<x_i$ for $x_i\in(0,0.5)$. When $b_i<1$, one has that $p(b_i,x_i)<x_i$ for $x_i\in(0,5,1)$ and $p(b_i,x_i)>x_i$ for $x_i\in(0,0.5)$.

{\color{black} We summarize the above observations in the following
\begin{enumerate}
    \item If $b_i>0$, then $p(b_i,x_i)>p(0,x_i)$ if $x_i\in(0.5,1)$, and $p(b_i,x_i)<p(0,x_i)$ if $x_i\in(0,0.5)$.
    \item If  $b_i = 1$, then $p(1,x_i)=x_i$.
    \item If $b_i > 1$, then $p(b_i,x_i)>x_i$ for $x_i\in(0.5,1)$ and $p(b_i,x_i)<x_i$ for $x_i\in(0,0.5)$.
    \item If $b_i < 1$, then $p(b_i,x_i)<x_i$ for $x_i\in(0.5,1)$ and $p(b_i,x_i)>x_i$ for $x_i\in(0,0.5)$.
\end{enumerate}
Item 1) indicates that individual $i$'s next opinion $x_i(k+1)$ under the bias model update rule \eqref{general} is closer to the polarized value of 0 (if $x_i(k) \in (0, 0.5)$) or 1 (if $x_i(k) \in (0.5, 1)$) when compared to $x_i(k+1) = p(0, x_i(k))$ of an individual $i$ described by the DeGroot model. It is by this mechanism that \eqref{general} captures an individual who, for $b_i > 0$, assimilates a balanced mixture of influence with a bias, more readily accepting neighboring information that supports his or her current opinion, while placing a lower weight on neighboring information that opposes his or her current opinion.

Item 2) illustrates a biased individual whose non-neutral opinion remains unchanged in the presence of equal information from both ends of the opinion spectrum; there is a perfect balance between biased assimilation and social influence from neighbors' opinions. Item 3) indicates that when $b_i > 1$, biased assimilation
overpowers the social influence, and the individual tends to an extreme opinion, \textit{even though the overall social influence due to the neighbors' opinions is unchanged}. This represents the scenario in which ``biased assimilation causes individuals to arrive at more extreme opinions after
being exposed to identical, inconclusive evidence" \cite{LoRoLe79}. Item 4) shows an individual whose the level of biased assimilation is not sufficient to overcome the mix of information from neighbors' opinions (social influence). Thus, $x_i$ tends to $0.5$, which is when the social influence from both ends of the spectrum is equal. Based on the above discussion, we say individual $i$ has weak bias if $b_i < 1$, or strong bias if $b_i > 1$, or intermediate bias if $b_i = 1$.}

\begin{remark}
\color{black} For some models, each individual has a parameter describing her susceptibility to external influence (the parameter is constant in \cite{FrJo99} and opinion-dependent in \cite{AmBuSi17}). However, both models share the same property; when an individual $i$ is exposed to two equal pieces of opinions from either end of the spectrum, the opinion furthest from opinion $x_i$ is more attractive. This contrasts our conclusion above; for an individual $i$ with $b_i > 0$, the opinion closer to opinion $x_i$ is more attractive.
\end{remark}


\section{Main Results} \label{result}

In this section, we present the theoretical results, interweaved with discussion and interpretation in the social context, with the proofs presented in the next section.
We will study the equilibria (and also their stability) of the system \eqref{general} for both $b_i>0$ and $b_i<0$. It turns out that this is a  challenging problem in general and some results we obtain only establish local stability.


Let
$$f_i(x)\triangleq\frac{w_{ii}x_i+(x_i)^{b_i}s_i}{w_{ii}+(x_i)^{b_i}s_i+(1-x_i)^{b_i}(d_i-s_i)}.$$
The update of the opinions of all $N$ individuals in the network is rewritten as
\begin{equation}\label{eq:system3}
x(k+1)= F(x(k)),
\end{equation}
where $x = [x_1, \hdots, x_N]^\top$ and $F(x) = [f_1(x), \hdots, f_N(x)]^\top$.

For system (\ref{general}) with $b_i>0$, note that if  $x_i(k) = 1$  ($x_i(k) = 0$), then $x_i(k^\prime) = 1$ ($x_i(k^\prime) = 0$) for all $k^\prime \geq k$.
It can be verified that for  $b_i>0$, we have that $\bm{0}_N,\ \bm{1}_N,$ and $\frac{1}{2}\bm{1}_N$ are equilibria of system (\ref{eq:system3}). We refer to $x^\ast = \bm{0}_N$ and $x^\ast = \bm{1}_N$ as extreme consensus and $x^\ast = \frac{1}{2}\bm{1}_N$ as neutral consensus.
Any vector with all entries either 0 or 1 is also an equilibrium; without loss of generality, we denote such an equilibrium as $[\bm{0}_{n_1}^T,\bm{1}_{n_2}^T]^T$ with $n_1+n_2=N$ and represents polarization of the network. Let the extreme, neutral consensus, and polarization equilibria of the system (\ref{eq:system3}) be respectively denoted by
\begin{eqnarray*}
&x_a^\ast=\bm{0}_N,\ &x_c^\ast=\bm{1}_N,\\
&x_d^\ast=\frac{1}{2}\bm{1}_N,\ &x_e^\ast=[\bm{0}_{n_1}^T,\bm{1}_{n_2}^T]^T.
\end{eqnarray*}

Besides the above equilibria, there may exist other equilibria of the system depending on the graph $\mathbb G$, the value of the bias parameter $b_i$, and the weights $w_{ij}$. We give examples in the sequel. If the bias parameter $b_i<0$ but is close to 0, then $x_d^\ast=\frac{1}{2}\bm{1}_N$ is an equilibrium of system (\ref{eq:system3}). Though a rigorous proof is missing, we conjecture the following based on numerous simulations.

{\color{black}{\bf Conjecture:} For a given network topology and initial states, if the system \eqref{eq:system3} with $b_i = 0$ for all $i \in \mathcal{V}$ converges (DeGroot model), then the system~(\ref{eq:system3}) with $b_i > 0$ for all $i \in \mathcal{V}$ will also converge.}

\subsection{Extreme and Neutral Consensus Equilibria}

We first discuss stability of equilibria corresponding to extreme consensus and neutral consensus.

\begin{theorem}\label{thm:1}
Suppose that the neighbor graph $\mathbb G$ is strongly connected {\color{black}and $b_i>0, \forall i \in \mathcal{V}$}. Then, $x_a^\ast=\bm{0}_N$ and $x_c^\ast=\bm{1}_N$ are locally exponentially stable equilibria and $x_d^\ast=\frac{1}{2}\bm{1}_N$ is an unstable equilibrium of system \rep{eq:system3}.
 \end{theorem}

In the social context, the result of Theorem~\ref{thm:1} indicates that individuals' biased assimilation makes it possible for a network \textit{to reach a
consensus that is more extreme} ($x(\infty) = \bm{1}_N$ and $x(\infty) = \bm{0}_N$) than any individual's initial opinion $x_i(0)$. For example,
one could have $x_i(0)\in(1-\epsilon_1,1-\epsilon_2)$ for all $i$, with sufficiently small $\epsilon_1>\epsilon_2>0$, and we get $x(\infty)= \bm{1}_N$, which means that $\max_i x_i(0) < \max_i x_i(\infty)$. This points to the dangers of biased assimilation in a network of individuals who all begin with similar opinions. One could say the network of individuals is
``self-extremizing". {\color{black}Theorem~\ref{thm:1} also tells us that when individuals exhibit biased assimilation, it is unlikely for a network to reach the unstable state of neutral consensus (in which every individual adopts the neutral opinion). However, it might be possible that stable equilibria exist in which a subset of the individuals (but not all) adopt the neutral opinion. For many established models \cite{De74,HeKr02,FrJo99,Al13}, the example initial states above will yield $\max_i x_i(0) \geq \max_i x_i(\infty)$. Some models \cite{MaFlKi14,Du17,ShPrJoBaJo17} can have $\max_i x_i(0) < \max_i x_i(\infty)$, but only exhibit extreme polarization (see Section~\ref{ssec:polarization_equib} below) exists in \cite{MaFlKi14,Du17,ShPrJoBaJo17}, and not extreme consensus.}

The paper \cite{DaGoLe13} showed that the biased assimilation model exhibits  extreme polarization $[\bm{0}_{n_1}^T,\bm{1}_{n_2}^T]^T$ in a two-island network, which is also not present in the existing models. {\color{black}Since \cite{DaGoLe13} requires that $n_1 \neq 0$ and $n_2 \neq 0$, this means that \cite{DaGoLe13} does not study the stability of extreme consensus states, as in our paper (extreme consensus can be considered as a special case of the polarization equilibria, with $n_1 = 0$ or $n_2 = 0$).}  We now detail a result on when the neutral consensus equilibrium is stable.

\begin{theorem}\label{thm:2}
Suppose that the neighbor graph $\mathbb G$ is strongly connected, $w_{ii}>0$ and {\color{black}$b_i \in [-\epsilon, 0)$, $\forall i \in \mathcal{V}$.} If $\epsilon > 0$ is sufficiently small, then $x_d^\ast=\frac{1}{2}\bm{1}_N$ is a locally exponentially stable equilibrium of the system \rep{eq:system3}.
\end{theorem}


To get some idea on the region of attraction of equilibria corresponding to extreme consensus (i.e. those reported in Theorem~\ref{thm:1}), we present the following result with the system starting from some restricted initial states.

\begin{theorem}\label{thm:convergence}
Consider the system \eqref{eq:system3}, and {\color{black}let $b_i \geq 1$ for all $i \in \mathcal{V}$}. Suppose that the neighbor graph $\mathbb G$ is strongly connected. Then,
\begin{enumerate}
\item
If $x_i(0)\ge 0.5$ for all $i \in \mathcal{V}$ and there exists at least one $j \in \mathcal{V}$
such that $x_j(0)>0.5$, then $x_i(k)$ will asymptotically converge to $1$ for all $i \in \mathcal{V}$.
\item
If $x_i(0)\le 0.5$ for all $i \in \mathcal{V}$ and there exists at least one $j\in \mathcal{V}$
such that $x_j(0)<0.5$, then $x_i(k)$ will asymptotically converge to $0$ for all $i \in \mathcal{V}$.
\end{enumerate}
\label{ji}\end{theorem}

From Theorem~\ref{thm:convergence}, we conclude that for $b_i \geq 1$, the region of attraction for extreme consensus is in fact quite large. In particular, for individuals with \textit{intermediate or strong levels of biased assimilation}, a network will ``self-extremize'' to a state of extreme consensus if all individuals begin on the same side of the opinion spectrum ($x_i(0) \geq 0.5$ or $x_i(0) \leq 0.5$ for all $i$), even if initially the individuals have varying degrees of support for the position at 1 or 0.  {\color{black}An echo chamber \cite{BaJoNaTuBo15} is a scenario whereby an individual only has access to information that supports his or her current opinion (this access may be a deliberate result of the individual's actions, or an unintended consequence of enabling technology, e.g. recommender systems). Theorem~\ref{thm:convergence} illustrates the dangerous consequence, viz. extreme consensus, of having individuals with intermediate/strong bias assimilation together in an echo chamber. }

\subsection{Polarization Equilibria}\label{ssec:polarization_equib}

For the stability of the equilibrium  $x_e^\ast=[\bm{0}_{n_1}^T,\bm{1}_{n_2}^T]^T$ corresponding to polarization, we have the following results on strongly connected graphs, complete graphs and two-island networks introduced below.


\begin{theorem}\label{thm:SCgraph}
If the neighbor graph $\mathbb{G}$ is strongly connected, and {\color{black}$b_i\in (0,1)$ for all $i \in \mathcal{V}$}, then the equilibrium $x_e^\ast = [\bm{0}_{n_1}^T,\bm{1}_{n_2}^T]^T$ of the system \eqref{eq:system3}  is unstable for every $n_1 = 1, \hdots, N-1$.
\end{theorem}

\begin{theorem}\label{thm:complete}
For an undirected complete neighbor graph $\mathbb{G}$ with weights $w_{ij}=1$, $i\neq j$, $i,j\in\mathcal V$, the equilibrium $x_e^\ast = [\bm{0}_{n_1}^T,\bm{1}_{n_2}^T]^T$ of the system \eqref{eq:system3}, with $n_1 = 2, \hdots, N-2$, {\color{black}is unstable when $b_i=1$ for all $i \in \mathcal{V}$ and is locally exponentially stable when $b_i>1$ for all $i \in \mathcal{V}$.}
\end{theorem}

Next we introduce the two-island network model studied in \cite{DaGoLe13}, which is used to model a {\em homophilous} network. Consider an undirected network in which the nodes in $\mathcal V$ are partitioned into two types, say $\tau_1,\tau_2.$  Let $\mathcal V_i$ denote the set of nodes of type $\tau_i$ and $|\mathcal V_i|=n_i,\ i=1,2$. Without loss of generality, assume that $\mathcal V_1=\{1,\ldots,n_1\}$ and $\mathcal V_2=\{n_1+1,\ldots,N\}.$ Assume that  each node in $\mathcal V_1$ has $n_1p_s$ neighbors in $\mathcal V_1$ and $n_1p_d$ neighbors in $\mathcal V_2$, and  each node in $\mathcal V_2$ has $n_2p_s$ neighbors in $\mathcal V_2$ and $n_2p_d$ neighbors in $\mathcal V_1$, where $p_s,p_d\in(0,1)$ and $n_1p_s,n_1p_d,n_2p_s,n_2p_d$ are all integers.

\begin{theorem}\label{thm:twoisland}
Suppose that the neighbor graph $\mathbb G$ is a connected two-island network and $w_{ij}=1$, $i\neq j$, $i,j\in\mathcal V$. Then, $x_e^\ast = [\bm{0}_{n_1}^T,\bm{1}_{n_2}^T]^T$
is a locally exponentially stable equilibrium of the system \eqref{eq:system3} when {\color{black}$b_i \geq1$ for all $i \in \mathcal{V}$. }
\end{theorem}

The above theorem results can be summarized in context as follows. Theorem~\ref{thm:SCgraph} establishes a result of particular interest when considered in conjunction with  Theorem~\ref{thm:1}. The results show that a network of individuals with weak bias, $b_i \in (0,1)$, can converge to an extreme consensus (which is undesirable), the same weak bias ensures that polarization (a different type of undesirable equilibrium) is an unstable phenomenon. Theorem~\ref{thm:SCgraph} also tells us it is unlikely for a network to converge to a polarized state if individuals are only weakly biased. The phenomenon of polarization is stable only when individuals have intermediate or strong levels of bias (Theorems~\ref{thm:complete} and \ref{thm:twoisland}).  Efforts to reduce polarization could therefore first focus on reducing individual bias as opposed to e.g. changing network structure or introducing agents into the network strategically.


\begin{remark}
The two-island network was also analyzed in \cite{DaGoLe13}, but with convergence to results secured for initial states restricted to satisfy $x_i(0) = x_0 \in (0.5, 1)$ if $i \in \mathcal{V}_1$ and $x_j(0) = 1- x_0$ if $j\in \mathcal{V}_2$. For convergence to $x_e^\ast = [\bm{0}_{n_1}^T,\bm{1}_{n_2}^T]^T$, \cite{DaGoLe13} also provided a relaxation on the initial states, requiring that $b_i=b \geq 1, \forall i\in\mathcal V$  and  $x_i(0) \geq 0.5 + \epsilon$ if $i \in \mathcal{V}_1$ and $x_j(0) \leq 0.5 - \epsilon$ if $j\in \mathcal{V}_2$, with $\epsilon$ being dependent on $p_s/p_d$. In contrast, we analyze the local stability and instability of the polarization equilibrium for varying values of the bias parameter $b_i$.  Based on numerous simulations where we sampled $x_i(0)$ from a uniform distribution in $[0,1]$, we observed that polarization occurs for a large set  of strongly connected network topologies, such as regular graphs, complete graphs, random graphs, and small-world graphs,  if $b_i$ is much larger than 1,   while for specific network topologies like path graphs and star graphs, polarization does not occur.
\end{remark}

The system (\ref{eq:system3}) can have other equilibria and can exhibit rich asymptotic behaviors as will be illustrated in Section~\ref{simulation}. The following theorem establishes a case when other types of equilibria of the system (\ref{eq:system3}) exist and their stability is discussed.

\begin{theorem}\label{thm:star}
Let $b_i=1, \forall\,i \in \mathcal{V}$. Consider an undirected star graph with the weights $w_{ij}=1$, $i\neq j$, $i,j\in\mathcal V$. Without loss of generality, suppose that node 1 is the center node. The equilibria of system (\ref{eq:system3}) 
include those vectors whose elements are either zero or one, and $x^\ast = [\frac{1}{2},a_2,\dots,a_N]^T$ with $a_i\in[0,1]$ and $\sum_{i=2}^Na_i=\frac{N-1}{2}$. If $N$ is odd, the system has  additional equilibria of the form $x^\ast = [c,\bm{0}_{\frac{N-1}{2}}^T,\bm{1}_{\frac{N-1}{2}}^T]^T$ with $c\in(0,1)$.
Among these equilibria, $x_a^\ast=\bm{0}_N,$ and $x_c^\ast=\bm{1}_N$  are locally exponentially stable and all the other equilibria are unstable.
\end{theorem}

Consider equilibria of the form
$x^\ast = [\frac{1}{2},a_2,\dots,a_N]^T$ with $\sum_{i=2}^Na_i= (N-1)/2$, or $x^\ast = [c,\bm{0}_{\frac{N-1}{2}}^T,\bm{1}_{\frac{N-1}{2}}^T]^T$ with $c\in(0,1)$. Theorem~\ref{thm:star} establishes that it is possible under biased assimilation to split a star network so that the leaf nodes separate into 2 nonempty factions, one supporting the opinion represented by 1, and the other supporting the opinion represented by 0. In fact, the support can be of varying levels of intensity, with different faction sizes, since one only requires that $a_i\in[0,1]$ and $\sum_{i=2}^Na_i=\frac{N-1}{2}$. The centre node acts as a ¡°mediating¡± individual to the two factions. However, such an equilibrium is unstable.

\section{Analyses}\label{analysis}


In this section, we prove the theorems in the previous section.
We linearize the system (\ref{eq:system3}) to analyze the local stability of these equilibria.
Let
$g_i(x)\triangleq w_{ii}+x_i^{b_i}s_i+(1-x_i)^{b_i}(d_i-s_i)$
for $i=1,\dots,N$. By calculation, one obtains that the Jacobian of $F(x(k))$ in (\ref{eq:system3}), $\frac{\partial F}{\partial x} = (\frac{\partial f_i}{\partial x_j})_{N\times N}$, has entries
\begin{align}\label{eq:jacob_diag_entry}
\frac{\partial f_i}{\partial x_i}=&\frac{1}{g_i^2(x)}
\bigg[(w_{ii}+b_i x_i^{b_i-1}s_i)g_i(x)-(w_{ii}x_i+x_i^{b_i} s_i)\nonumber \\
&\quad\times\left[b_i x_i^{b_i-1}s_i-b_i(1-x_i)^{b_i-1}(d_i-s_i)\right]\bigg]
\end{align}
and
\begin{align}\label{eq:jacob_offdiag_entry}
\frac{\partial f_i}{\partial x_l}=&\frac{1}{g_i^2(x)}
\bigg[x_i^{b_i}w_{il}g_i(x)\nonumber \\
&-(w_{ii}x_i+x_i^{b_i}s_i)(x_i^{b_i} w_{il}-(1-x_i)^{b_i} w_{il})\bigg]
\end{align}
for $l\neq  i$ and $i,l \in \mathcal{V}$.

A real matrix $M = (m_{ij})_{N\times N}$ is called a \emph{nonnegative} matrix  if $m_{ij}\geq0,\ i,j=1,\ldots,N$. The spectral radius of square matrix $M$ is denoted as $\rho(M)$. Before proving the theorems, the following lemma is first introduced that will be used later.
\begin{lemma}\label{lm:1}
(\cite{HoJo85})
Suppose $M\in \R^{N\times N}$ and $M$ is a nonnegative matrix. Then $\rho(M)$ is an eigenvalue of $M$ and $\min_{1\leq i\leq N}\sum_{j=1}^Nm_{ij}\leq \rho(M)\leq\max_{1\leq
i\leq N}\sum_{j=1}^Nm_{ij}$.
\end{lemma}


{\it Proof of Theorem~\ref{thm:1}:}
Consider the equilibrium $x_a^\ast=\bm{0}_N$. For all $i\in \mathcal{V}$, one knows that $s_i=0$ and
$g_i(x_a^\ast)=w_{ii}+d_i$.
One can derive using (\ref{eq:jacob_diag_entry}) and (\ref{eq:jacob_offdiag_entry}) that
\begin{equation*}
\frac{\partial f_i}{\partial x_i}\bigg|_{x_a^\ast}=\frac{w_{ii}}{g_i(x_a^\ast)}\,,\ \text{and }
\frac{\partial f_i}{\partial x_l}\bigg|_{x_a^\ast}=0,\ \text{for } l\neq  i.
\end{equation*}
Thus, the Jacobian matrix at the equilibrium $x_a^\ast=\bm{0}_N$ becomes
$$P\triangleq\frac{\partial F}{\partial x}\bigg|_{x_a^\ast}=\diag\left\{\frac{w_{11}}{g_1(x_a^\ast)},\frac{w_{22}}{g_2(x_a^\ast)},\ldots,
\frac{w_{NN}}{g_N(x_a^\ast)}\right\}.$$
Note that $g_i(x_a^\ast)=w_{ii}+\sum_{j\in\mathcal{N}_i}w_{ij}.$
The eigenvalues of $P$ are $w_{ii}/(w_{ii}+\sum_{j\in\mathcal{N}_i}w_{ij}),\ i\in\mathcal V$, which lie in the interval $[0,1)$ as long as each agent has at least one neighbor. Since $\mathbb G$ is strongly connected, $\rho(P)<1$ and thus the equilibrium $x_a^\ast=\bm{0}_N$ is locally exponentially stable.


For the equilibrium $x_c^\ast=\bm{1}_N$, observe that for $i\in\mathcal V$, one has
$g_i(x^\ast_c)=w_{ii}+d_i$. This yields
\begin{equation*}
\frac{\partial f_i}{\partial x_i}\bigg|_{x_c^\ast}=\frac{w_{ii}}{g_i(x_c^\ast)}\;,\text{ and }
\frac{\partial f_i}{\partial x_l}\bigg|_{x_c^\ast}=0,\, \text{for } l\neq  i
\end{equation*}
Thus the Jacobian matrix at the equilibrium $x_c^\ast=\bm{1}_N$ becomes
$$P\triangleq\frac{\partial F}{\partial x}\bigg|_{x_c^\ast}=\diag\left\{\frac{w_{11}}{g_1(x_c^\ast)},\frac{w_{22}}{g_2(x_c^\ast)},\ldots,
\frac{w_{NN}}{g_N(x_c^\ast)}\right\}.$$
The eigenvalues of $P$ are $w_{ii}/(w_{ii}+d_i),\ i\in\mathcal V$, which lie in the interval $[0,1)$ as in the previous case. Thus the equilibrium $x_c^\ast=\bm{1}_N$ is locally exponentially stable.


For the equilibrium $x_d^\ast=\frac{1}{2}\bm{1}_N$, one obtains $g_i(x_d^\ast)=w_{ii}+d_i/2^{b_i},$
for all $i\in\mathcal V$. It can be further calculated that
\begin{equation}\label{eq:jacob1}
\frac{\partial f_i}{\partial x_i}\bigg|_{x_d^\ast}=\frac{w_{ii}+\frac{b_id_i}{2^{b_i}}}{g_i(x_d^\ast)},\quad \frac{\partial f_i}{\partial x_l}\bigg|_{x_d^\ast}=\frac{w_{il}}{2^{b_i}g_i(x_d^\ast)},\;\text{for } l \neq i.
\end{equation}

The above implies that the Jacobian matrix $P\triangleq\frac{\partial F}{\partial x}|_{x_d^\ast}$ at $x_d^\ast=\frac{1}{2}\bm{1}_N$ is a nonnegative matrix.
Using Lemma \ref{lm:1} and \eqref{eq:jacob1}, one can compute that the spectral radius obeys
\begin{equation*}
\rho(P)\geq\min_{i=1,\ldots,N}\sum_{j=1}^Np_{ij}=1+\min_{i=1,\ldots,N}\frac{b_i\frac{1}{2^{b_i}}d_i}{w_{ii}+\frac{1}{2^{b_i}}d_i}>1.
\end{equation*}
Thus $x_d^\ast=\frac{1}{2}\bm{1}_N$ is an unstable equilibrium.
\hfill $\Box$


{\it Proof of Theorem~\ref{thm:2}:}
Similar calculations to the proof of Theorem~\ref{thm:1} shows that the Jacobian matrix evaluated at $x_d^\ast=\frac{1}{2}\bm{1}_N$, denoted $P \triangleq\frac{\partial F}{\partial x}|_{x_d^\ast}$, has the same entries as in (\ref{eq:jacob1}). The off-diagonal elements of $P$ are nonnegative. Since $w_{ii}>0$ and $b_i\in[-\epsilon,0)$, for all $i\in\mathcal V$, with $\epsilon > 0$ sufficiently small, one has
$$\frac{w_{ii}+b_i\frac{1}{2^{b_i}}d_i}{g_i(x_d^\ast)}
=\frac{w_{ii}+b_i\frac{1}{2^{b_i}}d_i}{w_{ii}+\frac{1}{2^{b_i}}d_i}\geq0$$
for all $i\in\mathcal V$, and hence $P$ is a nonnegative matrix. By Lemma~\ref{lm:1}, the spectral radius of $P$ satisfies that
$$\rho(P)\leq\max_{i=1,\ldots,N}\sum_{j=1}^Np_{ij}=1+\max_{i=1,\ldots,N}\frac{\frac{b_i}{2^{b_i}}\sum_{j\in\mathcal N_i}w_{ij}}{g_i(x_d^\ast)}<1.$$
Therefore $x_d^\ast=\frac{1}{2}\bm{1}_N$ is a locally exponentially stable equilibrium of the system (\ref{eq:system3}) for $b_i\in[-\epsilon,0)$ when $\epsilon$ is sufficiently small.
\hfill $\Box$




{\it Proof of Theorem~\ref{ji}:}
We first prove item 1). Consider any $i\in \mathcal V$, and observe that
\begin{align*}
&x_i(k+1) - x_i(k) \\
&= \frac{\zeta_i(x(k))}
{w_{ii}+(x_i(k))^{b_i}s_i(k)+(1-x_i(k))^{b_i}(d_i-s_i(k))}.
\end{align*}
where $\zeta_i(x) =  x_i^{b_i}s_i-x_i^{b_i+1} s_i-x_i(1-x_i)^{b_i}(d_i-s_i)$.

Proving that $x_i(k+1) - x_i(k) \geq 0$ is equivalent to proving that $\zeta_i(x(k)) \geq 0$ since the denominator of the equation above is positive. Rearranging the terms in $\zeta_i(x)$, and recalling that $d_i = \sum_{j\in \mathcal{N}_i} w_{ij}$ and $s_i = \sum_{j\in \mathcal{N}_i} w_{ij}x_j$, yields
\begin{align*}
    \zeta_i(x) & = \sum_{j\in\mathcal{N}_i} w_{ij} x_i^{b_i}\Big(x_j(1-x_i) - \frac{x_i}{x_i^{b_i}}(1-x_i)^{b_i}(1-x_j)\Big).
\end{align*}
Since $x_i \in [0.5,1],i\in\mathcal V$, $x_i^{b_i} > 0$, implying that $\zeta_i(x) \geq 0$ if $x_j(1-x_i)  - {x_i^{-b_i}}x_i(1-x_i)^{b_i}(1-x_j)  \geq 0$,
or equivalently:
\begin{align}
\frac{1- x_i}{x_i} \geq \left(\frac{1-x_i}{x_i}\right)^{b_i} \frac{1-x_j}{x_j} \label{eq:ben_01}
\end{align}
holds for all $j\in \mathcal{N}_i$. Trivially, \eqref{eq:ben_01} holds if $x_i = 1$, so let us consider $x_i \in [0.5, 1)$.  Notice that $x_i \in [0.5, 1) \Rightarrow (1-x_i)/x_i \leq 1$ with equality if and only if $x_i = 0.5$. Thus, \eqref{eq:ben_01} holds if $x_i \in [0.5, 1)$, with equality if and only if $x_j = 0.5$ and either (i) $b_i=1$ or (ii) $x_i = 0.5$. With $x_j \in [0.5, 1]$, $j\in \mathcal{N}_i$, we can then conclude that $\zeta_i(x) > 0$ if (i) $\exists j \in \mathcal{N}_i : x_j > 0.5$, or (ii) $b_i > 1$ and $x_i \in (0.5, 1)$. If $x_i(0)\ge 0.5$ for all $i\in\mathcal V$, then $x_i(k+1) \ge x_i(k)$ for all $i\in\mathcal V$ and all time $k$.
Moreover, since there exists at least one $j\in\mathcal V$
such that $x_j(0)>0.5$ and $\mathbb G$ is strongly connected, unless $x_i(0)=1$ for all $i\in\mathcal V$, there exists a $p\in\mathcal V$ such that $p\neq j$ and $x_p(1) > x_p(0)\ge 0.5$.
Repeating this argument, one concludes that there exists a finite $\tau$ such that $x_i(k)>0.5$ for all $i\in\mathcal V$ and $k\ge\tau$.

Consider the Lyapunov function
$V(x(k)) = 1-\min_{i\in\mathcal V}x_i(k)$.
From \rep{general}, if $x_i(k)=1$, then $x_i(k+1)=1$, which implies that
if $x_i(0) = 1$, then $x_i(k)=1$ for all time $k$.
Thus, if $x_i(k)=1$ for all $i\in\mathcal V$ at some time $k$, then
$V(x(k))=0$ and $V(x(k+1))=0$.
Suppose that
there exists at least one agent $p$ such that $x_p(k)<1$ at a specific time $k$. Without loss of generality, assume $k\ge \tau$. From the preceding discussion, $x_p(k+1)>x_p(k)$, which implies that
$\min_{i\in\mathcal V}x_i(k+1)>\min_{i\in\mathcal V}x_i(k)$, and thus  $V(x(k+1))<V(x(k))$.
By Lyapunov's stability theorem for discrete-time autonomous systems \cite[Theorem 13.2]{HaCh08}, $\lim_{k\to\infty} x_i(k) = 1$ asymptotically for all $i\in\mathcal V$.

Item 2) can be proved using arguments similar to those in the proof of item 1), with the Lyapunov function $V(x(k)) = \max_{i\in\mathcal V} x_i(k)$.
\hfill $\Box$


Before proving Theorems~\ref{thm:SCgraph}-\ref{thm:twoisland}, we calculate the elements of the Jacobian matrix of the system (\ref{eq:system3}) at the equilibrium $x_e^\ast=[\bm{0}_{n_1}^T,\bm{1}_{n_2}^T]^T$.


Let $x_{ei}^\ast$ be the $i$th entry of $x_{e}^\ast$, $\mathcal N_i^{(0)}=\{j:j\in\mathcal N_i, x_{ej}^\ast=0\}$ and $\mathcal N_i^{(1)}=\{j:j\in\mathcal N_i, x_{ej}^\ast=1\}$ for $i\in\mathcal V$, and $$d_i^{(0)}=\sum_{j\in\mathcal N_i^{(0)}}w_{ij},\ d_i^{(1)}=\sum_{j\in\mathcal N_i^{(1)}}w_{ij}.$$
For agent $i$ that satisfies $x_{ei}^\ast=0$, it is easy to see that
$g_i(x_e^\ast)=w_{ii}+d_i^{(0)}.$
Calculations show that
\begin{eqnarray*}\label{eq:pol1}
\frac{\partial f_i}{\partial x_i}\bigg|_{x_e^\ast}\!=\!\frac{(w_{ii}+{x_{ei}^\ast}^{b_i-1}b_id_i^{(1)})g_i(x_e^\ast)-{x_{ei}^\ast}^{2b_i-1}b_i(d_i^{(1)})^2}{g_i^2(x_e^\ast)},
\end{eqnarray*}
and $\frac{\partial f_i}{\partial x_l}|_{x_e^\ast}=\frac{1}{g_i^2(x_e^\ast)}{x_{ei}^\ast}^{b_i}=0$ for $l\neq  i$.
For agent $i$ that satisfies $x_{ei}^\ast=1$, one has
$g_i(x_e^\ast)=w_{ii}+d_i^{(1)}$.
Eq.~(\ref{eq:jacob_diag_entry}) then yields
\begin{eqnarray}\label{eq:pol3}
\frac{\partial f_i}{\partial x_i}\bigg|_{x_e^\ast} &=&\frac{1}{g_i(x_e^\ast)} \left[w_{ii}+(1-{x_{ei}^\ast})^{b_i-1}b_id_i^{(0)}\right],
\end{eqnarray}
and for $l\neq  i$, \eqref{eq:jacob_offdiag_entry} evaluates to be
\begin{equation}\label{eq:pol4}
\frac{\partial f_i}{\partial x_l}\bigg|_{x_e^\ast}\!\!=\frac{1}{g_i^2(x_e^\ast)}[w_{il}g_i(x_e^\ast)-w_{il}(w_{ii}+d_i^{(1)})]=0.
\end{equation}



{\it Proof of Theorem~\ref{thm:SCgraph}:} Since the graph is strongly connected, for any equilibrium $x_e^\ast=[\bm{0}_{n_1}^T,\bm{1}_{n_2}^T]^T$ with $n_1+n_2=N$ and  a given $n_1 \in \{1, \hdots, N-1\}$, there exists an agent $l$ such that $x^\ast_{el}=1$ and $d_l^{(0)}>0$. When $0<b_i<1$ for all $i \in \mathcal{V}$, it follows from \eqref{eq:pol3} that
$$\frac{\partial f_l}{\partial x_l}\bigg|_{x_e^\ast}=\frac{1}{g_l(x_e^\ast)} \left[w_{ll}+(1-{x_{el}^\ast})^{b_i-1}b_id_l^{(0)}\right]=+\infty.$$
In view of \eqref{eq:pol4}, the Jacobian matrix at the equilibrium $x_e^\ast=[\bm{0}_{n_1}^T,\bm{1}_{n_2}^T]^T$  is a diagonal matrix with at least one element equal to $+\infty$. Therefore, the equilibrium $x_e^\ast=[\bm{0}_{n_1}^T,\bm{1}_{n_2}^T]^T$  is unstable. \hfill $\Box$


{\it Proof of Theorem~\ref{thm:complete}:}
When $b_i=1$ for all $i\in \mathcal{V}$, for agent $i$ with $x_{ei}^\ast=0$ and agent $l$ with $x_{el}^\ast=1$, one has
\begin{equation*}
\frac{\partial f_i}{\partial x_i}\bigg|_{x_e^\ast}=\frac{w_{ii}+d_i^{(1)}}{w_{ii}+d_i^{(0)}} \,,\text{\,and }\; \frac{\partial f_l}{\partial x_l}\bigg|_{x_e^\ast}=\frac{w_{ll}+d_l^{(0)}}{w_{ll}+d_l^{(1)}},
\end{equation*}
respectively. Then the Jacobian matrix $P=\frac{\partial F}{\partial x}|_{x_e^\ast}$ at the equilibrium $x_e^\ast=[\bm{0}_{n_1}^T,\bm{1}_{n_2}^T]^T$ is
\begin{align}
P&=\!\diag\!\left\{\!\frac{w_{11}\!+\!d_1^{(1)}}{w_{11}\!+\!d_1^{(0)}},
\frac{w_{22}\!+\!d_2^{(1)}}{w_{22}\!+\!d_2^{(0)}},\ldots,\frac{w_{NN}\!+\!d_N^{(0)}}{w_{NN}\!+\!d_N^{(1)}}\!\right\}.\label{eq:Pb=1}
\end{align}
Suppose that $n_1\geq n_2$. For the $i$-th agent with $x_{ei}^\ast=1$, in view of (\ref{eq:Pb=1}), the $i$-th diagonal element of $P$ is given by $p_{ii}=(w_{ii}+n_1)/(w_{ii}+n_2-1)>1$.
If $n_1< n_2$, one can similarly show that there exists a diagonal element of $P$ that is greater than 1. In both cases, $P$ has an eigenvalue greater than 1. Therefore the  equilibrium  $x_e^\ast=[\bm{0}_{n_1}^T,\bm{1}_{n_2}^T]^T$ is unstable when $b_i=1$ for all $i\in \mathcal{V}$.

When $b_i>1$ for all $i\in \mathcal{V}$, for agent $i$ with $x_{ei}^\ast=0$, one has $\frac{\partial f_i}{\partial x_i}|_{x_e^\ast}=w_{ii}/(w_{ii}+d_i^{(0)})$.
Since the graph is complete and $n_1\geq2$, $d_i^{(0)}>0$ and therefore $0\leq w_{ii}/(w_{ii}+d_i^{(0)})<1$. For agent $l$ with $x_{el}^\ast=1$, one has $\frac{\partial f_l}{\partial x_l}|_{x_e^\ast}=w_{ll}/(w_{ll}+d_l^{(1)})$.
Similarly, one derives that  $0\leq w_{ll}/(w_{ll}+d_l^{(1)})<1$.
Then the Jacobian matrix $P=\frac{\partial F}{\partial x}|_{x_e^\ast}$ evaluated at the equilibrium $x_e^\ast=[\bm{0}_{n_1}^T,\bm{1}_{n_2}^T]^T$ is
\begin{align}
P=\!\diag\!\left\{\!\frac{w_{11}}{w_{11}\!+\!d_1^{(0)}},\frac{w_{22}}{w_{22}\!+\!d_2^{(0)}},
\ldots,\frac{w_{NN}}{w_{NN}\!+\!d_N^{(1)}}\!\right\}\label{eq:Pblarge}
\end{align}
and has spectral radius $\rho(P)<1.$ It follows that $x_e^\ast=[\bm{0}_{n_1}^T,\bm{1}_{n_2}^T]^T$ is locally exponentially stable when $b_i>1$ for all $i\in \mathcal{V}$.
\hfill $\Box$



{\it Proof of Theorem~\ref{thm:twoisland}:}
From the definition of the two-island network model, the following inequalities hold \begin{equation}\label{eq:b=1}
d_1^{(0)}>d_1^{(1)},
d_2^{(0)}>d_2^{(1)},\ldots,
d_N^{(1)}>d_N^{(0)}.
\end{equation}
For $b_i=1$ and $b_i>1$, the Jacobian matries are given by \eqref{eq:Pb=1} and (\ref{eq:Pblarge}), respectively. In both cases, one can see that the eigenvalues of $P$ lie in the interval $[0,1)$ and thus $x_e^\ast=[\bm{0}_{n_1}^T,\bm{1}_{n_2}^T]^T$ is locally exponentially stable.
\hfill $\Box$

{\it Proof of Theorem~\ref{thm:star}:} Let $x^\ast=[x_1^\ast,x_2^\ast,\dots,x_N^\ast]^T$ be an equilibrium of the system (\ref{eq:system3}). If $x_1^\ast=1$, then one can show that $x_i^\ast, i=2,\ldots, N,$ can only be either 0 or 1. The same conclusion holds for the case when $x_1^\ast=0$. Hence $\bm{1}_N$, $\bm{0}_N$, and $[\bm{0}_{n_1}^T,\bm{1}_{n_2}^T]^T$ with $n_1+n_2=N,$ are equilibria of the system.

Suppose that $x_1^\ast\neq1$ and $x_1^\ast\neq0$. For the center node, it follows from
$$\frac{w_{11}x_1^\ast+x_1^\ast\sum_{i=2}^Nx_i^\ast}
{w_{11}+x_1^\ast\sum_{i=2}^Nx_i^\ast+(1-x_1^\ast)(N-1-\sum_{i=2}^Nx_i^\ast)}
=x_1^\ast$$
 that $\sum_{i=2}^Nx_i^\ast=\frac{N-1}{2}$. For the nodes $i=2,\dots,N$, it should hold that
\begin{equation}\label{eq:star1}
\frac{w_{ii}x_i^\ast+x_i^\ast x_1^\ast}
{w_{ii}+x_i^\ast x_1^\ast+(1-x_i^\ast)(1-x_1^\ast)}
=x_i^\ast.
\end{equation}
Suppose there exists some $i \in \{2, \hdots, N\}$ such that $x_i^\ast\neq0$ and $x_i^\ast\neq 1$. Then, (\ref{eq:star1})  holds if $x_1^\ast=\frac{1}{2}$. In conclusion, $[\frac{1}{2},a_2,\dots,a_N]^T$ with $\sum_{j=2}^Na_i=\frac{N-1}{2}$ and $a_i\geq0, i\in \{2, \hdots, N\}$, are equilibria of the system. Moreover, if $N$ is odd, the system has the additional equilibria $[c,\bm{0}_{\frac{N-1}{2}}^T,\bm{1}_{\frac{N-1}{2}}^T]^T$ with $c\in(0,1)$.

We first discuss the stability of the polarization equilibria. Note that node 1 is the center. Consider an equilibrium whose  first element is  0 and there exists some other element, say $i$, whose value is 1. Then according to the Jacobian matrix $P$ given by (\ref{eq:Pb=1}), its $i$-th diagonal element is $p_{ii}=(w_{ii}+d_i^{(0)})/(w_{ii}+d_i^{(1)})=(w_{ii}+1)/w_{ii}$.
Since $p_{ii}>1$, such an equilibrium is unstable. The stability of an equilibrium with the first element equal to 1 can be similarly discussed.

We now check the stability of the equilibria $x^\ast = [\frac{1}{2},a_2,\dots,a_N]^T$. The Jacobian matrix $P$ at the  equilibria $x^\ast =[\frac{1}{2},a_2,\dots,a_N]^T$ is given by
$$\begin{bmatrix}
    1 & \frac{1}{2w_{11}+N-1} & \dots & \frac{1}{2w_{11}+N-1} \\
    \frac{4a_2(1-a_2)}{2w_{22}+1} & 1 & \dots & 0 \\
    \vdots & \vdots & \ddots & \vdots \\
    \frac{4a_N(1-a_N)}{2w_{NN}+1} & 0 & \dots & 1 \\
  \end{bmatrix},$$
which  is nonnegative. Suppose without loss of generality that $0<a_i<1$, $2\leq i\leq k,$  for some $2\leq k\leq N$, and $a_j=0$ or $a_j=1$, for $k+1\leq j\leq N$. Then the leading principle submatrix $P_k$  of order $k$ of $P$  is irreducible and can be written as
$$P_k=
\begin{bmatrix}
    1 & \frac{1}{2w_{11}+N-1} & \dots & \frac{1}{2w_{11}+N-1} \\
    \frac{4a_2(1-a_2)}{2w_{22}+1} & 1 & \dots & 0 \\
    \vdots & \vdots & \ddots & \vdots \\
    \frac{4a_k(1-a_k)}{2w_{kk}+1} & 0 & \dots & 1 \\
  \end{bmatrix}.$$
From Lemma~\ref{lm:1}, it follows that $\rho(P_k)$ is an eigenvalue of $P_k$ and
\begin{equation*}
\rho(P_k)\geq\min\left\{\!1\!+\!\frac{k\!-\!1}{N\!-\!1},1\!+\!\min_{2\leq i\leq k}\{1+4a_i(1-a_i)\}\right\}\!>\!1.
\end{equation*}
Since $\rho(P_k)$ is an eigenvalue of $P$ as well, the spectral radius $\rho(P)>1$ and thus the equilibria $x^\ast =[\frac{1}{2},a_2,\dots,a_N]^T$ are unstable.

Now consider the case when $N$ is odd. The Jacobian $P$ at the equilibria $x^\ast =[c,\bm{0}_{\frac{N-1}{2}}^T,\bm{1}_{\frac{N-1}{2}}^T]^T$ is given by
\begin{equation*}
\begin{bmatrix}
    1 & \frac{4c(1-c)}{2w_{11}+N-1} & \dots & \frac{4c(1-c)}{2w_{11}+N-1} \\
    0 & \frac{w_{22}+c}{w_{22}+1-c} & \dots & 0 \\
    \vdots & \vdots & \ddots & \vdots \\
    0 & 0 & \dots & \frac{w_{NN}+1-c}{w_{NN}+c} \\
  \end{bmatrix},
\end{equation*}
which is a nonnegative matrix. If $c\neq \frac{1}{2}$, then either $\frac{w_{22}+c}{w_{22}+1-c}$ or $\frac{w_{NN}+1-c}{w_{NN}+c}$ will be larger than 1. Therefore the  equilibria $x^\ast =[c,\bm{0}_{\frac{N-1}{2}}^T,\bm{1}_{\frac{N-1}{2}}^T]^T$ with $c\neq \frac{1}{2}$ are unstable.

For  the equilibrium $x^\ast =[\frac{1}{2},\bm{0}_{\frac{N-1}{2}}^T,\bm{1}_{\frac{N-1}{2}}^T]^T$, consider a small perturbation around this equilibrium. Take the initial condition of the system (\ref{eq:system3}) as $x(0)=[c,a\bm{1}_{\frac{N-1}{2}}^T,\bm{1}_{\frac{N-1}{2}}^T]^T$, where $c>\frac{1}{2}$ is close to $\frac{1}{2}$ and $a>0$ is close to 0. It is clear that for the agents $i=(N+1)/{2}+1,\ldots,N$, $x_i(k)=1$ for all $k\geq0$. For the center node,
$$x_1(1)=\frac{w_{11}c+c(1+a)\frac{N-1}{2}}{w_{11}+c(1+a)\frac{N-1}{2}+(1-c)(1-a)\frac{N-1}{2}}.$$
One can show that  $x_1(1)>c$ as long as $\frac{1}{2}<c<1$ and $a>0$. For $i=2,\ldots,(N+1)/{2}$,
$$x_i(1)=\frac{w_{ii}a+ac}{w_{ii}+ac+(1-a)(1-c)},$$
and $x_i(1)>a$ as long as $0<a<1$ and $c>\frac{1}{2}$. By induction, the states of the agents $i=1,\ldots,(N+1)/{2}$ are strictly monotonically increasing and will  converge to 1 as time goes to infinity. The equilibrium $[\frac{1}{2},\bm{0}_{\frac{N-1}{2}}^T,\bm{1}_{\frac{N-1}{2}}^T]^T$ is unstable.
\hfill $\Box$


\section{Numerical simulations}\label{simulation}

In this section, we perform  several simulations to show the rich asymptotic behaviors of the system (\ref{eq:system3}), including some equilibria not studied in Section~\ref{result}. {\color{black}In each of the following simulations, a two-island network model with each island consisting of 50 nodes is considered. For each node, the number of neighbors of the same type is $n_1p_s = n_2 p_s = 4$ and the number of neighbors of  the other type is $n_1p_d = n_2 p_d =2$. Edges are bidirectional, i.e. $(j,i)\in \mathcal{E}$ then $(i,j)\in\mathcal{E}$, but the weights $w_{ij}$ and $w_{ji}$ are not necessarily equal, thus making the graph directed. In particular, if $(j,i) \in \mathcal{E}$, we drawn $w_{ij}$ from a uniform distribution with interval $[0.5,1.5]$, and set $w_{ii}=0$ for all $i\in\mathcal V$. }


{\color{black}In the first case, we consider when $b_i$ for all $i\in\mathcal V$
are chosen randomly from a uniform distribution in the interval $[1.01,1.5]$, i.e. all individuals have strong bias. The initial states of the agents are chosen randomly from a uniform distribution in the interval $[0,1]$, and the evolution of the states of the agents are illustrated in Fig.~\ref{fig:two_island101}, from which one can see that the system reaches an extreme polarization equilibrium. If $b_i$ for all $i\in\mathcal V$ are much larger than 1, we observe from extensive simulations that extreme polarization is also observed for a large class  of strongly connected network topologies  such as regular graphs, random graphs, and small-world graphs. This illustrates the important role of individuals with strong bias in creating a polarized network state.}

\begin{figure} [htbp]
\begin{center}
\includegraphics[width=8cm]{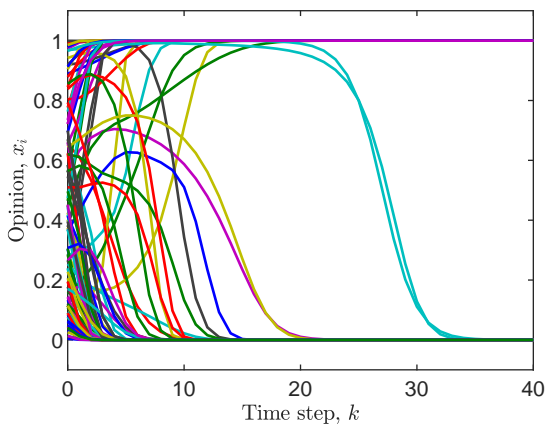}    
\caption{The system state under a two-island network  with $b\in[1.01,1.5]$ and randomized initial states.}\label{fig:two_island101}
\end{center}
\end{figure}

{\color{black}To illustrate that there are other equilibria which are very different to those analyzed in Section~\ref{result}, we present the following simulations. We now draw $b_i$ for all $\ i\in\mathcal V$ from a uniform distribution of interval $\in[0.5,1.5]$, so that some individuals have weak bias and some have strong bias. If the initial states are randomly chosen from a uniform distribution in $[0,1]$, then Fig.~\ref{fig:two_island1} illustrates that the states of most of the agents converge either to 0 or 1 and the final states of the remaining agents lie in between.  Again, similar results to Fig.~\ref{fig:two_island1} can be observed in other network topology types, including path networks, regular networks and small-world networks.}

\begin{figure} [htbp]
\begin{center}
\includegraphics[width=8cm]{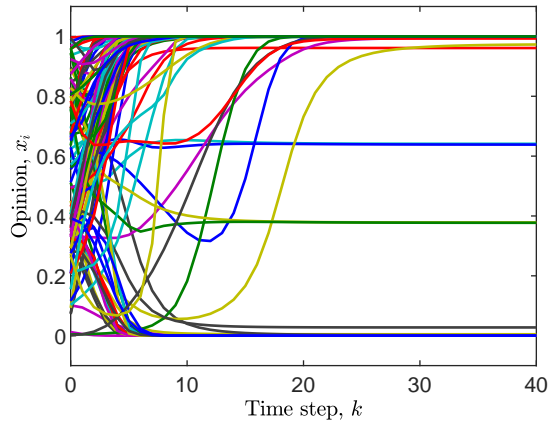}
\caption{The system state under a two-island network with $b_i\in[0.5,1.5]$ and randomized initial states.}\label{fig:two_island1}
\end{center}
\end{figure}

{\color{black}Now consider the case where $b_i$ for all $\ i\in\mathcal V$ belongs to a uniform distribution of interval $\in(0,1)$. For initial states uniformly randomly chosen from the interval $[0,1]$, two situations are typically observed for the state evolution of the system. In the first situation, the states of all agents converge either to 1 or 0, and in the other situation, the states of most of the agents converge to values close either to 0 or 1 and the final states of the remaining agents lie in $[0,1]$. As all $b_i$ values tend closer to 0, the situation in which the states of all the agents converge to an extreme consensus equilibrium occurs more frequently. When $b_i$ is close to 1, for some specific initial states, the agents converge to two clusters of opinions close to the extreme polarization equilibria. For example, we consider $b_i\in[0.8,0.9]$ under the two-island network, and the initial states of agents 1 to 50 are randomly chosen from a uniform distribution of interval $[0.15,0.2]$ and the remaining agents have initial states from a uniform distribution of interval $[0.75,0.8]$. Fig.~\ref{fig:bipartite} shows that the network converges to a steady state in which the two islands have states that are close to the extreme values of 0 and 1.}

\begin{remark}
\color{black}We have shown that there are equilibria other than those studied in Section~\ref{result}. Although not shown, we also observed that heterogeneous $b_i$ can generate equilibria that does not exist for a homogeneous bias parameter. Similarly, there may be equilibria for undirected networks which do not exist for directed graphs, and vice versa.
\end{remark}

\begin{figure} [htbp]
\begin{center}
\includegraphics[width=8cm]{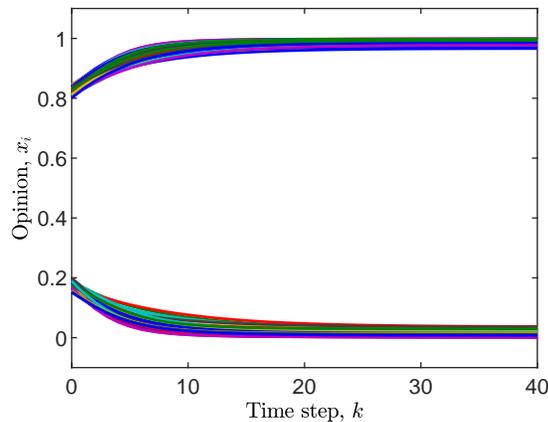}    
\caption{The system state under a two-island network with $b_i\in[0.8,0.9]$ and the initial state $x_i(0)\in [0.15,0.2],i=1,\dots,50$ and $x_i(0)\in [0.8,0.85],i=51,\dots,100$.}\label{fig:bipartite}
\end{center}
\end{figure}

\section{Conclusion} \label{conclusion}

In this paper, we have studied the equilibria and their stability for a recently proposed nonlinear opinion dynamics model with biased assimilation in which each agent is associated with a bias parameter. {\color{black}We have shown that, with heterogeneous bias parameter values, the stability of certain equilibria depend on the degree of bias and the topology of the neighbor relationships among the agents. Both theoretical analyses and numerical simulations have shown that both the value of the bias parameter and the network topology play a key role in determining the limiting opinion distributed in the network. For future work, we aim to further study the region of attraction of the different equilibria and explore the general convergence condition for arbitrary strongly connected networks and arbitrary initial states, though a conjecture was given in Section~\ref{result}.}


\bibliographystyle{plain}        
\bibliography{bias_arxivV2}           

\begin{thebibliography}{10}

\bibitem{Aj01}
I.~Ajzen.
\newblock Nature and operation of attitudes.
\newblock {\em Annual Review of Psychology}, 52:27--58, 2001.

\bibitem{Al13}
C.~Altafini.
\newblock Consensus problems on networks with antagonistic interactions.
\newblock {\em IEEE Transactions on Automatic Control}, 58(4):935--946, 2013.

\bibitem{AmBuSi17}
V.~Amelkin, F.~Bullo, and A.~K. Singh.
\newblock Polar opinion dynamics in social networks.
\newblock {\em IEEE Transactions on Automatic Control}, 62(11):5650--5665,
  2017.

\bibitem{AnYe19}
B.~D.~O. Anderson and M.~Ye.
\newblock Recent advances in the modelling and analysis of opinion dynamics on
  influence networks.
\newblock {\em International Journal of Automation and Computing},
  16(2):129--149, 2019.

\bibitem{An81}
N.~H. Anderson.
\newblock {\em Foundations of Information Integration Theory}.
\newblock Academic Press, 1981.

\bibitem{BaJoNaTuBo15}
P.~Barber{\'a}, J.~T. Jost, J.~Nagler, J.~A. Tucker, and R.~Bonneau.
\newblock Tweeting from left to right: {I}s online political communication more
  than an echo chamber?
\newblock {\em Psychological Science}, 26(10):1531--1542, 2015.

\bibitem{BlHeTs09}
V.~D. Blondel, J.~M. Hendrickx, and J.~N. Tsitsiklis.
\newblock On {K}rause's multi-agent consensus model with state-dependent
  connectivity.
\newblock {\em IEEE Transactions on Automatic Control}, 54(11):2586--2597,
  2009.

\bibitem{ChQiLiQiBuSh19}
Z.~Chen, J.~Qin, B.~Li, H.~Qi, P.~Buchhorn, and G.~Shi.
\newblock {Dynamics of Opinions with Social Biases}.
\newblock {\em Automatica}, 106(8):374--383, 2019.

\bibitem{DaGoLe13}
P.~Dandekar, A.~Goel, and D.~T. Lee.
\newblock Biased assimilation, homophily, and the dynamics of polarization.
\newblock {\em Proceedings of the National Academy of Sciences},
  110:5791--5796, 2013.

\bibitem{De74}
M.~H. DeGroot.
\newblock Reaching a consensus.
\newblock {\em Journal of the American Statistical Association},
  69(345):118--121, 1974.

\bibitem{Po14}
M.~Dimock, J.~Kiley, S.~Keeter, and C.~Doherty.
\newblock {Political Polarization in the American Public}.
\newblock Technical report, Pew Research Center, 2014.

\bibitem{Du17}
P.~Duggins.
\newblock A psychologically-motivated model of opinion change with applications
  to american politics.
\newblock {\em Journal of Artificial Societies and Social Simulation},
  20(1):1--13, 2017.

\bibitem{FrJo90}
N.~E. Friedkin and E.~C. Johnsen.
\newblock Social influence and opinions.
\newblock {\em Journal of Mathematical Sociology}, 15(3-4):193--205, 1990.

\bibitem{FrJo99}
N.~E. Friedkin and E.~C. Johnsen.
\newblock Social influence networks and opinion change.
\newblock {\em Advances in Group Processes}, 16(1):1--29, 1999.

\bibitem{ShPrJoBaJo17}
{G. Shi}, A.~Proutiere, M.~Johansson, J.~S. Baras, and K.~H. Johansson.
\newblock Emergent behaviors over signed random dynamical networks:
  {R}elative-state-flipping model.
\newblock {\em IEEE Transactions on Control of Network Systems}, 4(2):369--379,
  2017.

\bibitem{HaCh08}
W.~M. Haddad and V.~Chellaboina.
\newblock {\em Nonlinear Dynamical Systems and Control: A Lyapunov-Based
  Approach}.
\newblock Princeton University Press, Princeton, 2008.

\bibitem{HeKr02}
R.~Hegselmann and U.~Krause.
\newblock Opinion dynamics and bounded confidence: {M}odels, analysis and
  simulation.
\newblock {\em Journal of Artificial Societies and Social Simulation}, 5:1--24,
  2002.

\bibitem{HoJo85}
R.~A. Horn and C.~R. Johnson.
\newblock {\em Matrix Analysis}.
\newblock Cambridge University Press, Cambridge, U.K., 1985.

\bibitem{JiMiFrBu15}
P.~Jia, A.~MirTabatabaei, N.~E. Friedkin, and F.~Bullo.
\newblock Opinion dynamics and the evolution of social power in influence
  networks.
\newblock {\em SIAM Review}, 57(3):367--397, 2015.

\bibitem{LiChBaBe17}
J.~Liu, X.~Chen, T.~Ba\c{s}ar, and M.-A. Belabbas.
\newblock Exponential convergence of the discrete- and continuous-time
  {A}ltafini models.
\newblock {\em IEEE Transactions on Automatic Control}, 62(12):6168--6182,
  2017.

\bibitem{LoRoLe79}
C.~G. Lord, L.~Ross, and M.~R. Lepper.
\newblock Biased assimilation and attitude polarization: {T}he effects of prior
  theories on subsequently considered evidence.
\newblock {\em Journal of Personality and Social Psychology},
  37(11):2098--2109, 1979.

\bibitem{MaFlKi14}
M.~M\"{a}s, A.~Flache, and J.~A. Kitts.
\newblock Cultural integration and differentiation in groups and organizations.
\newblock In {\em Perspectives on Culture and Agent-based Simulations}, pages
  71--90. Springer, 2014.

\bibitem{MuDiLoFaGrPe02}
G.~D. Munro, P.~H. Ditto, L.~K. Lockhart, A.~Fagerlin, M.~Gready, and
  E.~Peterson.
\newblock Biased assimilation of sociopolitical arguments: Evaluating the 1996
  {U.S.} presidential debate.
\newblock {\em Basic and Applied Social Psychology}, 24(1):15--26, 2002.

\bibitem{PrRo17}
A.~V. Proskurnikov and R.~Tempo.
\newblock {A tutorial on modeling and analysis of dynamic social networks. Part
  I}.
\newblock {\em Annual Reviews in Control}, 43:65--79, 2017.

\bibitem{TaLo06}
C.~S. Taber and M.~Lodge.
\newblock {Motivated skepticism in the evaluation of political beliefs}.
\newblock {\em American Journal of Political Science}, 50(3):755--769, 2006.

\bibitem{XiCaJo16}
W.~Xia, M.~Cao, and K.~H. Johansson.
\newblock Structural balance and opinion separation in trust-mistrust social
  networks.
\newblock {\em IEEE Transcations on Control of Network Systems}, 3(1):46--56,
  2016.

\bibitem{YeLiAnYuBa19}
M.~Ye, J.~Liu, B.~D.~O. Anderson, C.~Yu, and T.~Ba{\c{s}}ar.
\newblock {Evolution of Social Power in Social Networks with Dynamic Topology}.
\newblock {\em IEEE Transaction on Automatic Control}, 63(11):3793--3808,
  November 2019.

\end{thebibliography}


\end{document}